\documentclass[a4paper,10pt]{amsart}

\usepackage{mescommandes}

\newtheorem{theo}{Theorem}
\newtheorem{pro}[theo]{Proposition}
\newtheorem{lem}[theo]{Lemma}
\newtheorem{cor}[theo]{Corollary}

\newcommand{\ddt}{\frac{d}{dt}}
\newcommand{\dif}[1]{\frac{d}{d#1}}
\newcommand{\df}[1]{\frac{d}{d#1}\Bigr|_{#1=0}}

\newcommand{\dpi}{\mathrm{d}\pi}
\newcommand{\odpi}{\ov{\mathrm{d}\pi}}

\title{Integrating representations of Banach--Lie algebras}
\author{St\'ephane Merigon}

\begin{document}
\maketitle

\section{Introduction}
In this note we give two integrability criteria for a
representation $\ha:\fg\ra\End(\D)$
of a Banach--Lie algebra $\fg$ as skew-symmetric unbounded 
operators on a dense domain $\D$ of a Hilbert space, that is, sufficient conditions 
for the existence of
a continuous unitary representation of a simply connected 
Banach--Lie group with Lie algebra $\fg$ (when it exists)
whose derived representation extends $\ha$.

After Nelson's famous criterion \cite{NEL}, new ones 
(for finite dimensional Lie algebras) appeared in the late sixties and early seventies, also based on analytic vectors (see \cite{FSSS,SI})
or on smooth vectors (see \cite{JOMO}). In both cases the key result is the validity of
the commutation relation
\begin{equation}\label{E:comrel}
e^{\ha(x)}\ha(y)e^{-\ha(x)}=\ha(e^{\ad x}y),
\end{equation}
for every $x$ in a Lie-generating subset of $\fg$ and any $y$ in 
$\fg$. It is used to transfer computations in the space of operators
to computations in the Lie algebra so that the integrability follows from formulas in $\fg$.

We consider a Banach--Lie algebra $\fg$ which decomposes as 
$$\fg=\fa_1\oplus\fa_2\oplus\dots\oplus\fa_n,$$
where $\fa_j$, $j=1,2,\dots n,$ are closed subspaces
and we prove that the representation integrates if \eqref{E:comrel} holds
for every $x\in\cup\fa_j$ and $y\in\fg$. This also leads to a criterion based on
analytic vectors and generalising \cite{FSSS}: the representation integrates if $\D$ consists
of analytic vectors for every $\ha(x)$, $x\in\cup \fa_j$
(for a generalisation of Nelson's criterion see \cite{LA}). The new difficulty 
in the Banach setting is that we have to differentiate
pathes of the form
\begin{equation}\label{E:newdef}
t\mt e^{\ha(x(t))}v,
\end{equation}
where $x(t)$ is a smooth path in $\fa_j$, but \eqref{E:comrel} 
enables us to do that. The derivative of \eqref{E:newdef} 
involves the logarithmic derivative of $x(t)$, therefore a formula
relating such derivatives is needed.

These results are particularly well-suited for the study of symmetric or
$3$-graded Lie algebras. In a forthcoming work 
\cite{MERNEEB}
with Karl-Hermann Neeb
we will use them to prove a generalisation of the
L\"uscher--Mack Theorem \cite[Appendix C]{LM} to the Banach--Lie setting. This
is a crucial tool for the study of the representation theory 
of the automorphisms groups of
infinite dimensional real symmetric domains. 

The integrability criteria are stated in Section~\ref{S:2}.
In Section~\ref{S:3} the relation~\eqref{E:comrel} is discussed and used to prove the differentiability of \eqref{E:newdef} while the relevant formula for logarithmic derivatives is given in 
Section~\ref{S:4}. The proof of the integrability criteria is 
achieved in Section \ref{S:5}.

\section{Main results}\label{S:2}

Let $\fg$ be a Banach--Lie algebra with Lie bracket $[\cdot,\cdot]$.
A representation $\ha$ of $\fg$ on a (dense) subspace $\D$ of a Hilbert space $\cH$
is a linear map which associates to any $x\in\fg$ a skew-symmetric unbounded operator $\ha(x):\D\ra\D$, in such a way that $$\ha([x,y])=[\ha(x),\ha(y)]:=\ha(x)\ha(y)-\ha(y)\ha(x).$$ 
The representation is said to be 
\emph{strongly continuous} if for every $v\in\D$ the map $$\ha^v:\fg\ra\cH,\ x\mt\ha(x)v$$ is
continuous. 

Let $G$ be a Banach--Lie group with Lie algebra $\fg$ and 
$\pi:G\ra\dU(\cH)$ be a continuous unitary
representation of $G$ in $\cH$. Each element $x\in\fg$ gives rise
to a one-parameter unitary group $\pi(\exp tx)$ and hence by Stone's
theorem to a skew-adjoint operator
$$\odpi(x)v:=\df{t}\pi(\exp tx)v$$
defined on the set of vectors $v$ for which the limit exists.
The derived representation is
the (strongly continuous) representation of $\fg$ defined on the space $\cH^\infty$ of smooth
vectors (those vectors for which the orbit map is smooth) by 
$$\mathrm{d}\pi(x):=\odpi(x)|_{\cH^\infty}.$$
When $G$ is finite dimensional $\cH^\infty$ is dense and since it is
invariant under the action of $G$, the operators $\dpi(x)$ are 
essentially skew-adjoint, but for an arbitrary Banach-Lie group $\cH^\infty$ may be empty. For issues related to smooth vectors for representations of infinite dimensional Lie groups
see \cite{NEE2}.

We say that a strongly continuous representation $\ha$ on the
dense domain $\D$ \emph{integrates} to a continuous unitary representation of $G$ if there exists such a representation $\pi$  with $\D\subseteq\cH^\infty$ and $\mathrm{d}\pi|_{\D}=\ha$.

We below assume that $\fg$ has a 
decomposition
$$\fg=\fa_1\oplus\fa_2\oplus\dots\oplus\fa_n,$$
where $\fa_j$, $j=1,2,\dots n,$ are closed subspaces.

The main theorem will be stated for a strongly continuous representation $\ha$ on a dense domain $\D$ which
satisfies the following assumptions:
\begin{enumerate}
 \item[{\bf (A1)}] For all 
$x\in\cup\fa_j$, $\ha(x)$ is essentially 
skew-adjoint, i.e., its closure $\ov{\ha(x)}$ is skew-adjoint, hence generates
a strongly continuous one-parameter unitary group 
$e^{t\ha(x)}:=e^{t\ov{\ha(x)}}$, $t\in\bR$.
\item[{\bf (A2)}]  For all $x\in\cup\fa_j$, 
$e^{\ha(x)}\D\subseteq \D$.
 \item[{\bf (A3)}]  For all $(x,y)\in 
(\cup\fa_j,\fg)$ and $v\in\D$, the commutation relation
$$e^{\ha(x)}\ha(y)e^{-\ha(x)}v=\ha(e^{\ad x}y)v$$
holds.
\end{enumerate}

\begin{theo}\label{T} 
Let $G$ be a simply
connected Banach--Lie group with Lie algebra $\fg$. Any strongly 
continuous representation of $\fg$ satisfying {\bf (A1-3)} 
integrates to a continuous unitary representation of 
$G$.
\end{theo}

First we give a corollary in which the assumption {\bf (A3)} is weakened 
(see \cite[Theorem 9.1]{JOMO},
where in the case of finite dimensional Lie algebras it is weakened even further).

\begin{cor}\label{C:C1} 
Let $G$ be a simply
connected Banach--Lie group with Lie algebra $\fg$. Any strongly 
continuous representation of $\fg$ satisfying {\bf (A1-2)} and
 \begin{enumerate}
 \item[{\bf (A3')}]  For all $(x,y)\in 
(\cup\fa_j,\fg)$ and $v\in\D$, the map
$\bR\ra\D,\ t \mt \ha(y)e^{t\ha(x)}v$
is continuous.
\end{enumerate} 
integrates to a continuous unitary representation of 
$G$.
\end{cor}

The second corollary is an integrability criterion based on analytic vectors. 
It generalises the Integrability Theorem \cite[6.8]{NEE2}, 
where $\D$ is assumed to consists of analytic vectors
for every $\ha(x)$, $x\in\fg$, but the techniques involved are
completely different. For finite dimensional Lie algebras it was proved by
M. Flato, J. Simon, H. Snellman, and  D. Sternheimer (see \cite{FSSS}). Note that J. Simon also
proved
\cite{SI} that is it sufficent to assume that $\D$ consists of analytic vectors for every $\ha(x)$, $x\in S$, where $S$ is a Lie-generating subset of $\fg$.

\begin{cor}\label{C:C2} Let $G$ be a simply
connected Banach--Lie group with Lie algebra $\fg$. Let $\ha$ be a strongly 
continuous representation of $\fg$ 
over a dense domain which consists of analytic vectors for
the operators $\ha(x)$, $x\in\cup\fa_j$. Then $\ha$ integrates to a continuous unitary representation of $G$.
\end{cor}

In the next section we show how the corollaries follow from Theorem~\ref{T},
which will be proved in the last section.

\section{The commutation relation}\label{S:3}

In this section we first prove that the assumptions of
Corollary~\ref{C:C1}, as well as those of Corollary~\ref{C:C2},
lead to the assumptions {\bf (A1-3)}. 
Those results can be found for finite dimensional Lie algebras in \cite[Ch. 3]{JOMO} and in \cite{FSSS} respectively, but although
they extend directly to the Banach setting we give full proofs for the sake of completeness. Then we show that the commutation relation \eqref{E:comrel}
implies the differentiability of a family of operators which is crucial in the proof of the main theorem.


We will need the following product rule.

\begin{lem}\label{L:PR} Let $E$ and $F$ be two Banach spaces and
$L_s(E,F)$ denote the space of continuous linear operators from $E$ to $F$ endowed
with the strong operator topology.
Let $t\mt K(t)\in L_s(E,F)$ be a continuous
path such
that $t\mt K(t)v$ is differentiable for every $v$ in a subspace $\D$ of $E$ and let $\hg(t)$ be a differentiable
path in $\D$. We write $K'(t):\D\ra F$ for the linear
operator obtained by $K'(t)v:=\ddt K(t)v$, $v\in\D$.
Then $t\mt K(t)\hg(t)$ is differentiable with
$$\ddt K(t)\hg(t)=K'(t)\hg(t)+K(t)\hg'(t).$$
\end{lem}
\begin{proof} We write
\begin{multline*}
\frac{1}{h}(K(t+h)\hg(t+h)-K(t)\hg(t))=\\
\frac{1}{h}\big(K(t+h)\hg(t+h)-K(t+h)\hg(t)\big)+\frac{1}{h}
\big(K(t+h)\hg(t)-K(t)\hg(t)\big).
\end{multline*}
The second term converges to $K'(t)\hg(t)$ and the first one is 
equal to
$$K(t+h)\left(\frac{\hg(t+h)-\hg(t)}h-\hg'(t)\right)+K(t+h)\hg
'(t).$$
which converges to $K(t)\hg'(t)$ by the Principle of Uniform Boundedness: 
If $\C$ is a compact neighbourhood of $t$, then for every $v\in\cH$, $\sup_{s\in\C}K(s)v$ is bounded and hence $\sup_{s\in\C}\nm{K(s)}$ is bounded.
\end{proof}

\begin{lem}\label{L:formder} Consider two unbounded operators $A$ and $B$ defined on a dense domain
$\D$ of the Hilbert space $\cH$. Assume that $A$ is essentially skew-adjoint, that $A\D\subseteq\D$ and
$e^{tA}\D\subseteq\D$, and that $B$ is closable. Let $v\in\D$ such that
$t\mt BAe^{tA}v$ is continuous. Then $t\mt Be^{tA}v$ is differentiable with
$$\dif{t}Be^{tA}v=BAe^{tA}v.$$
\end{lem}
\begin{proof}
We have
$$e^{tA}v-v=\int_0^t{Ae^{sA}vds}.$$
Let $\ov B$ be the closure of $B$. Its domain $\D(\ov{B})$ is a Banach space
when endowed with the graph norm $\nm{w}_B:=\nm{w}+\nm{\ov{B}w}$, where $\nm{\cdot}$ is the Hilbert norm in $\cH$,
and then 
$\ov{B}:\D(\ov{B})\ra\cH$ is a continuous linear operator.
By assumption the map $s\mt Ae^{sA}v$ is continuous for the graph norm and hence
the integral $\int_0^t{Ae^{sA}vds}$ exists 
as a Riemann integral in $\D(\ov{B})$. We therefore have
$$Be^{tA}v-Bv=\int_0^t{BAe^{sA}vds},$$
and the claim follows.
\end{proof}

\begin{pro}\label{P:com_rel} Let $\ha$ be a strongly 
continuous representation of the	
Banach--Lie algebra $\fg$ over a dense domain $\D$ of a Hilbert space
$\cH$. Let $x\in\fg$ such that
$\ha(x)$ is essentially skew-adjoint and such that the 
associated one-parameter unitary group leaves $\D$ invariant.
Let $\fa$ be a closed subspace
of $\fg$ which invariant under $\ad x$. Assume that for every
$v\in\D$ and every $y\in\fa$
the map $\bR\ra\cH,\ t\mt\ha(y)e^{t\ha(x)}v$ is continuous. 
Then we have for every $y\in\fa$ the commutation 
relation
$$e^{\ha(x)}\ha(y)e^{-\ha(x)}=\ha(e^{\ad x}y).$$
\end{pro}
\begin{proof} (See \cite[3.2 and 3.3]{JOMO})
Let $v\in\D$. We want to prove 
that the map
$$[0,1]\ra\cH,\ s\mt e^{(1-s)\ha(x)}\ha(e^{s\ad 
x}y)e^{(s-1)\ha(x)}v$$ 
is constant.
We will apply Lemma~\ref{L:PR} with
$$K(s):\fa\ra\cH,\ z\mt 
e^{(1-s)\ha(x)}\ha(z)e^{(s-1)\ha(x)}v$$
and $\hg(s)=e^{s\ad x}y$.
The continuity of the operator $K(s)$ 
follows directly from the strong continuity of $\ha$.
The continuity of the map $s\mt\ha(z)e^{(s-1)\ha(x)}\ha(x)v$
implies by Lemma~\ref{L:formder} the differentiability of 
$s\mt\ha(z)e^{(s-1)\ha(x)}v$, and the derivative is
$$\dif{s}\ha(z)e^{(s-1)\ha(x)}v=\ha(z)\ha(x)e^{(s-1)\ha(x)}v.$$
Hence, by Lemma~\ref{L:PR},
the map $s\mt 
K(s)z=e^{(1-s)\ha(x)}\ha(z)e^{(s-1)\ha(x)}v$ 
is differentiable with derivative
$$K'(s)z=e^{(1-s)\ha(x)}[\ha(z),\ha(x)]e^{(s-1)\ha(x)}v
.$$
Since
$s\mt e^{s\ad x}y$ is analytic and $\dif{s}e^{s\ad x}y=[x,e^{s\ad x}y]$, again by Lemma~\ref{L:PR},
we have 
\begin{multline*}
\dif{s}e^{(1-s)\ha(x)}\ha(e^{s\ad x}y)e^{(s-1)\ha(x)}v=\\
e^{(1-s)\ha(x)}[\ha(e^{s\ad x}y),\ha(x)]e^{(s-1)\ha(x)}v+ 
e^{(1-s)\ha(x)}\ha([x,e^{s\ad x}y])e^{(s-1)\ha(x)}v=0.
\end{multline*}
\end{proof}


Let $\ha:\fg\ra\End(\D)$ be a strongly continuous representation of
$\fg$ by skew-symmetric operators and let
$x\in\fg$ such that $\D$ consists 
of analytic vectors for $\ha(x)$. By Nelson's Theorem \cite[Lemma 5.1]{NEL} $\ha(x)$ then is essentially skew-adjoint. The key result here is \cite[Proposition 1]{FSSS},
the proof of which carries over to the Banach case without any change.

\begin{pro}[\cite{FSSS}]\label{P:FSSS} Let $\ha$ and $\ha'$
be strongly continuous representations of the Banach--Lie algebra
$\fg$ over the domains
$\D$ and $\D'$ (respectively), dense in $\cH$, with $\D\subseteq\D'$,
and such that for any $y\in\fg$, $\ha(y)$ is the restriction to
$\D$ of $\ha'(y)$. Then, if $\D$ is a domain of analytic vectors for
some $\ha(x)$, $x\in\fg$, we have for any $v\in\D$, $w\in\D'$, denoting by $\langle\cdot,\cdot\rangle$ the scalar product in $\cH$,
\begin{equation*}
\langle-e^{\ha(x)}\ha(y)v,w\rangle=\langle 
e^{\ha(x)}v,\ha'(e^{\ad x}y)w\rangle.
\end{equation*}
\end{pro}
Remark that if for every $y\in\fg$ one is given on $\D'\supseteq\D$ an unbounded skew-symmetric operator $\ha'(y):\D'\ra\D'$, so that
$\ha'(y)|_\D=\ha(y)$, then $\ha'$ is automatically a strongly continuous Lie algebra
homomorphism. Indeed, for $v\in\D$ and $w\in\D'$, we have
$$\langle-\ha(y)v,w\rangle=\langle 
v,\ha'(y)w\rangle.$$
Hence
\begin{multline*}
\langle v,\ha'(\lambda x+y)w\rangle=\langle 
-\ha(\lambda x+y)v,w\rangle=\\
\langle -\lambda\ha(x)v-\ha(y)v,w\rangle
=\langle v,\lambda\ha'(x)w+\ha'(y)w\rangle
\end{multline*}
and
\begin{multline*}
\langle v,\ha'([x,y])w\rangle=\langle 
\ha([y,x])v,w\rangle=\\
\langle \ha(y)\ha(x)-\ha(x)\ha(y)v,w\rangle
=\langle v,\ha'(x)\ha'(y)-\ha'(y)\ha'(x)w\rangle
\end{multline*}
show that $\ha'$ is a Lie algebra homomorphism.
Moreover a representation of Banach--Lie algebra is strongly continuous if and only if it is weakly continuous (\cite[Lemma~4.2]{NEE}). The domain $\D$ being dense, the set of functionals $\{\langle\cdot,v\rangle,\ v\in\D\}$
separates the points in $\cH$ and hence $\ha'$ is strongly continuous.

Assume now that $\ha(y)$ is essentially skew-adjoint for every $y\in\cup\fa_j$
and let
\begin{equation}\label{E:D'}
\D':=\bigcap_{\ell\in\bN,\, y_k\in\bigcup\fa_j}\D(\ov{\ha(y_\ell)}\dots\ov{\ha(y_1)}),
\end{equation}
where $\ov{\ha(y)}$ denotes the closure of $\ha(y)$.
Then $\D'$ contains $\D$ and is invariant under $\ov{\ha(y)}$, 
$y\in\cup\fa_j$, 
so we can set for such $y$,
$$\ha'(y):=\ov{\ha(y)}|_{\D'},$$
and for $y=y_1+\dots+y_n\in\fg$, $y_j\in\fa_j$, 
$$\ha'(y)=\ha'(y_1)+\dots+\ha'(y_n).$$
By the preceding remark this defines a strongly continuous representation of
$\fg$ extending $\ha$.
Applying Proposition~\ref{P:FSSS} several times, we see, since $\ha(y)^*=
\ov{\ha(y)}$ for $y\in\cup\fa_j$, that
for every $y_1,\dots,y_\ell\in\cup\fa_j$ and $w\in\D'$, $e^{\ha(x)}w\in\D(\ov{\ha(y_\ell)}\dots\ov{\ha(y_1)})$,
and hence
$$e^{\ha(x)}\D'\subseteq\D'.$$
Then Proposition~\ref{P:FSSS} also shows that the commutation relation holds on $\D'$:
\begin{cor}\label{C:comrel}
Let $\ha$ be a strongly 
continuous representation of the	
Banach--Lie algebra $\fg$ over a dense domain $\D$ such that $\ha(y)$
is essentially skew-adjoint for every $y\in\cup\fa_j$. Let $x\in\fg$ such that $\D$ consists of analytic vectors for
$\ha(x)$. Let $\D'\supseteq\D$ be the domain defined in \eqref{E:D'}
and $\ha'$ the corresponding extension of $\ha$.
Then $e^{\ha(x)}$
leaves $\D'$ invariant
and we have for every $y\in\fg$ the commutation 
relation
$$e^{\ha(x)}\ha'(y)e^{-\ha(x)}=\ha'(e^{\ad x}y).$$
\end{cor}


For the next proposition we will need the following lemma:

\begin{lem}\label{L:diff}
Consider two essentially skew-adjoint operators $A$ and $B$
defined on a common dense domain $\D$ 
and assume that for all $s\in\bR$, $e^{sA}\D\subseteq\D$. Let $v\in\D$ be such that 
$s\mt Be^{sA}v$
is continuous. Then
$$e^{tB}v-e^{tA}v=\int_0^t{e^{sB}(B-A)e^{(t-s)A}vds}.$$
\end{lem}
\begin{proof}
We apply Lemma~\ref{L:PR}  with $K(s)=e^{sB}$ and 
$\hg(s)=e^{(t-s)A}v$ to obtain
$$\dif{s}e^{sB}e^{(t-s)A}v=e^{sB}(B-A)e^{(t-s)A}v.$$
By assumption the right hand side is continuous. Hence  
the claim follows from the
Fundamental Theorem of Calculus.
\end{proof}
\begin{pro}\label{P:derpath} Let $\ha$ be a strongly 
continuous representation of the	
Banach--Lie algebra $\fg$ on a dense domain $\D$.
Let $I$ be a real interval and $I\ra\fg$, $t~\mt~x(t)$ be a continuous path such that
each $\ha(x(t))$ is essentially skew-adjoint and 
$e^{s\ha(x(t))}\D\subseteq\D$, $s\in\bR$. If, for every $t$, we have for sufficiently small $h$ and every $s\in \bR$
the commutation relation
\begin{equation}\label{E:localcomrel}
e^{s\ha(x(t))}\ha(x(t+h))e^{-s\ha(x(t))}=\ha(e^{s\ad x(t)}x(t+h)), \end{equation}
then for 
every $v\in\D$ the map 
$(s,t)\mt e^{s\ha(x(t))}v$ is continuous. If moreover the path $x(t)$ is differentiable, then $t\mt e^{\ha(x(t))}v$ is differentiabe with
$$\ddt 
e^{\ha(x(t))}v=\ha\left(\int_0^1{e^{s\ad{x(t)}}x'(t)ds}\right)e^{\ha(x(
t))}v.$$
\end{pro}
\begin{proof} Let us fix $t\in I$ and let $t+\N_t$ be a convex
neighbourhood of $t$ in $I$ such that
for $h\in\N_t$ the relation~\eqref{E:localcomrel} holds.
Let $v\in\D$. Rewriting \eqref{E:localcomrel} as
$$\ha(x(t+h))e^{s\ha(x(t))}v
=e^{s\ha(x(t))}\ha(e^{-s\ad{x(t)}}x(t+h))v,$$
we see that $s\mt\ha(x(t+h))e^{s\ha(x(t))}v$ is continuous.
We can therefore apply Lemma~\ref{L:diff} to obtain 
\begin{align*}
&e^{s\ha(x(t+h))}v-e^{s\ha(x(t))}v=
\int_0^s{e^{u\ha(x(t+h))}\ha\left(x(t+h)-x(t)\right)e^{(s-u)\ha(x(t))}vdu}\\
&=\int_0^s{e^{u\ha(x(t+h))}e^{(s-u)\ha(x(t))}\ha\left(e^{(u-s)\ad x(t)}(x(t+h)-x(t))\right)e^{(s-u)\ha(x(t))}vdu}.
\end{align*}
Thus, writing $\nm{\cdot}$ for the norm in $\fg$,
$$\hnm{e^{s\ha(x(t+h))}v-e^{s\ha(x(t))}v}\leq\hnm{s}\hnm{\ha^v}
e^{\hnm{s}\nm{\ad x(t)}}\nm{x(t+h)-x(t)}$$
and $(s,t)\mt e^{s\ha(x(t))}v$ is continuous.
Now assume that $t\mt x(t)$ is differentiable. Let us write
$$\frac{e^{\ha(x(t+h))}v-e^{\ha(x(t))}v}{h}=
\int_0^1{e^{s\ha(x(t+h))}\ha\left(\frac{x(t+h)-x(t)}{h}\right)e^{(1-s)\ha(x(t))}vds}$$
and let us define on $\N_t$ the function
\begin{equation*}
z(h)=
\begin{cases}
\frac{x(t+h)-x(t)}{h} & \text{for}\ h\neq0,\\
x'(t) & \text{for}\ h=0.
\end{cases}
\end{equation*}
The formula
$$\ha(z(h))e^{(1-s)\ha(x(t))}v=e^{(1-s)\ha(x(t))}\ha(e^{(s-1)
\ad x(t)}z(h))v,$$
which holds for $h\neq0$ by assumption and for $h=0$ by continuity,
shows that
$$(s,h)\mt \ha(z(h))e^{(1-s)\ha(x(t))}v$$ 
is continuous, and hence
$$(s,h)\mt e^{s\ha(x(t+h))}\ha(z(h))e^{(1-s)\ha(x(t))}v$$
is continuous. We can therefore pass to the 
limit under the integral sign to derive that
$$\ddt 
e^{\ha(x(t))}v=\int_0^1{e^{s\ha(x(t))}\ha(x'(t))e^{(1-s)\ha(x(
t))}vds}=\int_0^1{\ha(e^{s\ad{x(t)}}x'(t))e^{\ha(x(
t))}vds},$$
and the claim follows from the linearity and the continuity of $\ha$.
\end{proof}

\section{The right-logarithmic derivative}\label{S:4}

Let $\fg$ be a Banach--Lie algebra. Let $\U\subset\fg$ be a symmetric 
starlike neighbourhood of $0$ in $\fg$ such that the Dynkin series $x*y$ 
converges in $\U\times\U$. Then, for $x\in\U$, the maps
$$\lambda_xy:=x*y,\quad\hr_xy:=y*x\quad\text{and}\quad c_xy:=x*y*(-x)$$
are local diffeormorphisms at the origin. Moreover the differential at $0$
of $c_x$ 
is given by
$$Dc_x(0)=e^{\ad x}.$$

Let $I$ denote an interval of the real line. The \emph{right logarithmic derivative} (see \cite[II.4]{NEE3}) of a smooth path $\hg:I\ra\U$ is defined
by
$$\hd(\hg)_t=D\hr_{\hg(t)}(0)^{-1}\hg'(t).$$

\begin{lem} Let $x\in\U$ and $\hg(t)=tx$. Then
$$\hd(\hg)_t=x.$$
\end{lem}
\begin{proof}
We have $D\hr_{tx}(0)x=\lim_{h\ra0}\frac{hx*tx-tx}h=\lim_{h\ra0}\frac{(h+t)x-tx}h=x$.
\end{proof}

\begin{lem} \label{L:PRD} Let $\ha,\hb:I\mt\U$ two differentiable paths such that
$(\ha*\hb)(t):=\ha(t)*\hb(t)\in\U$. Then
$$\hd(\ha*\hb)_t=\hd(\ha)_t+e^{\ad \ha(t)}\hd(\hb)_t.$$
\end{lem}
\begin{proof} We have
\begin{align*}
\hd(\ha*\hb)_t&=D\hr_{(\ha*\hb)(t)}(0)^{-1}(\ha*\hb)'(t)\\
&=D\hr_{\ha(t)}(0)^{-1}D\hr_{\hb(t)}(\ha(t))^{-1}
\big(D\hr_{\hb(t)}(\ha(t))\ha'(t)+D\hl_{\ha(t)}(\hb(t))
\hb'(t)\big)\\
&=D\hr_{\ha(t)}(0)^{-1}\ha'(t)+D\hr_{\ha(t)}(0)^{-1}
D\hr_{\hb(t)}(\ha(t))^{-1}
D\hl_{\ha(t)}(\hb(t))\hb'(t)\\
&=D\hr_{\ha(t)}(0)^{-1}\ha'(t)+Dc_{\ha(t)}(0)
D\hr_{\hb(t)}(0)^{-1}\hb'(t)
\\
&=\hd(\ha)_t+e^{\ad \ha(t)}\hd(\hb)_t.
\end{align*}
\end{proof}
\noindent The next lemma says that the logarithmic derivative is
the pull-back of the Maurer-Cartan form on $\fg$ (and may therefore
be defined for any path in $\fg$). 
\begin{lem}
Let $\ha:I\mt\U$ be a differentiable path. Then
$$\hd(\ha)_t=\int_0^1{e^{s\ad{\ha(t)}}\ha'(t)ds}.$$
\end{lem}
\begin{proof}
Let us fix $t\in I$ and consider the map 
$$\psi:[0,1]\ra\fg,\quad s\mt\hd(s\ha)_t.$$
Then 
$$\psi(s+h)=\hd((s+h)\ha)_t=\hd((s\ha)*(h\ha))_t,$$
so we obtain with Lemma~\ref{L:PRD},
$$\psi(s+h)=\hd(s\ha)_t+e^{s\ad{\ha(t)}}\hd(h\ha)_t.$$
Hence we have
$$\psi'(s)=\lim_{h\ra 0}\frac{\psi(s+h)-\psi(s)}{h}=\lim_{h\ra 0}e^{s\ad{\ha(t)}}\frac1h\hd(h\ha)_t=e^{s\ad{\ha(t)}}\ha'(t),$$
and the result follows by integration.
\end{proof}

Assume now that $\fg$ decomposes as 
$$\fg=\fa_1\oplus\fa_2\oplus\dots\oplus\fa_n,$$
where $\fa_j$, $j=1,2,\dots n$ are closed subspaces. Then, for
every $j=1,2,\dots n$, there exists a $0$-neighbourhood $\V_j$ in $\fa_j$ such that
the map 
$$\V_1\times\V_2\times\dots\times\V_n\ra\fg,\quad (x_1,x_2,\dots,x_n)\mt x_1*x_2*\dots*x_n$$
is a diffeomorphism onto is image.
From now on we assume that we have chosen $\U$ starlike and small enough so that it is contained
in this image. So if $x,y\in\U$ then 
$$(tx)*y=x_1(t)*x_2(t)*\dots*x_n(t)$$
where $t\mt x_j(t)\in\fa_j$ is analytic.
\begin{pro}\label{P:strucform}
We have
$$x=\sum_{j=1}^{n}e^{\ad x_1(t)}\dots e^{\ad x_{j-1}(t)} 
\int_0^1{e^{s\ad x_j(t)}x_j'(t)ds}$$
\end{pro}
\begin{proof}
Let $\hg(t)=(tx)*y=x_1(t)*x_2(t)*\dots*x_n(t)$. 
The result follows by computing, using the preceding lemmas, the right logarithmic derivative of $\hg$ in its two expressions.
\end{proof}

\section{Proof of the main theorem}\label{S:5}

Let us consider a strongly continuous representation $\ha$ of $\fg$ on a dense domain $\D$ of the Hilbert space $\cH$ which
satisfies the assumptions {\bf (A1-3)} and
recall the
notations of the 
preceding section.

If 
$$z=z_1*\dots*z_n\in\V_1*\dots*\V_n$$
we set
\begin{equation}\label{E:constr}
\pi(z):=e^{\ha(z_1)}e^{\ha(z_2)}\dots e^{\ha(z_n)}.
\end{equation}
Let $\U'$ be a starlike $0$-neighbourhood in $\fg$ so that $\U'*\U'\subseteq\U$.
Then it suffices to show that for every $x,y\in\U'$,
\begin{equation}\label{F:mult}
\pi(x*y)=\pi(x)\pi(y),
\end{equation}
see, e.g., \cite[Ch. 3, \S 6, Lemma 1.1]{BOUR}.

Let us write
$$(tx)*y=x_1(t)*x_2(t)*\dots*x_n(t),\quad x_j(t)\in\fa_j,$$
so that
$$\pi((tx)*y)=e^{\ha(x_1(t))}e^{\ha(x_2(t))}\dots e^{\ha(x_n(t))}.$$
Let $v\in\D$ and 
$$\hg(t)=\pi((tx)*y)v.$$
Thanks to
Proposition~\ref{P:derpath}, we can use
Lemma~\ref{L:PR} several times to see that the map $t\mt\hg(t)$ is differentiable with
$$\hg'(t)=\sum_{j=1}^n{e^{\ha(x_1(t))}\dots 
e^{\ha(x_{j-1}(t))}
\ha\left(\int_0^1{e^{s\ad{x_j(t)}}x_j'(t)ds}\right)e^{\ha(x_{j}(t))}\dots 
e^{\ha(x_n(t))}v}.$$
Then repeated use of the commutation relation yields
$$\hg'(t)=\ha\left(\sum_{j=1}^ne^{\ad{x_1(t)}}\dots e^{\ad{x_{j-1}(t)}}\int_0^1{e^{s\ad{x_j(t)}}x_j'(t)}\right)
\hg(t),$$
which, according to Proposition~\ref{P:strucform}, amounts to
\begin{equation}\label{E:der}
\hg'(t)=\ha(x)\hg(t).
\end{equation}
But is it well known (cf. \cite[P. 431]{FSSS}) that the solution of the initial value problem
$$\hg(t)\in\D,\quad \hg'(t)=\ha(x)\hg(t), \quad \ha(0)=\pi(y)v,$$
is unique (and given by $e^{t\ha(x)}\pi(y)v$). Therefore we have
$$\pi(tx)\pi(y)v=\pi((tx)*y)v,$$ 
and this equality extends to $\cH$ to give
\eqref{F:mult} by evaluation at $t=1$. 
We also derive from~\eqref{E:der}, with $y=0$, that
$\ov{\dpi}(x)|_\D=\ha(x)$, and since by construction
$\pi(tx)\D\subseteq\D$, 
the operator $\ha(x)$ is essentially skew-adjoint, 
i.e. $\ov{\ha(x)}=\ov{\dpi}(x)$.
It also follows from \eqref{E:constr} that \eqref{E:comrel} holds for
every $x\in\fg$. Hence for $v\in\D$ the map
$$\fg\ra\cH,\ x\mt e^{\ov{\ha(x)}}v$$
is continuous (see Proposition~\ref{P:derpath}), and this implies that the representation $\pi$ is continuous.
Now we have 
$$\D\subseteq\D_\fg^\infty:=\bigcap_{\ell\in\bN,\, y_k\in\fg}\D(\ov{\dpi}(y_\ell)\dots\ov{\dpi}(y_1)),$$
but we know by \cite[Lemma 3.4, Remark 8.3]{NEE2} that $\D_\fg^\infty$ coincides with
the space of smooth vectors for $\pi$. 
This conludes the proof.

\section*{Acknowledgement} I thank Karl-Hermann Neeb for guiding me
through the literature
and for reading preliminary versions of the manuscript.

\bibliographystyle{amsalpha}
\bibliography{ref}
\end{document}